\documentclass[12pt]{article}

\makeatletter
\renewcommand{\@oddfoot}{\hfill \thepage}
\makeatother

\usepackage[T2A]{fontenc}
\usepackage[cp1251]{inputenc}
\usepackage[russian, english]{babel}
\usepackage{mathrsfs}
\usepackage{amssymb}
\usepackage{amsmath,amsthm,amsfonts}
\usepackage{cite}
\usepackage[dvipdf]{graphicx}
\usepackage{fancyhdr}
\usepackage[titletoc]{appendix}

\begin{document}

\begin{center}
{\bf ASYMPTOTIC RELATION FOR THE TRANSITION DENSITY \\ OF THE THREE-DIMENSIONAL MARKOV RANDOM FLIGHT\\ 
ON SMALL TIME INTERVALS}
\end{center}

\begin{center}
Alexander D. KOLESNIK\\
Institute of Mathematics and Computer Science\\
Academy Street 5, Kishinev 2028, Moldova\\
E-Mail: kolesnik@math.md
\end{center}

\vskip 0.2cm

\begin{abstract}
We consider the Markov random flight $\bold X(t), \; t>0,$ in the three-dimensional Euclidean space $\Bbb R^3$ with constant finite speed $c>0$ and the uniform choice of the initial and each new direction at random time instants that form a homogeneous Poisson flow 
of rate $\lambda>0$. Series representations for the conditional characteristic functions of $\bold X(t)$ corresponding to 
two and three changes of direction, are obtained. Based on these results, an asymptotic formula, as $t\to 0$, for the unconditional characteristic function of $\bold X(t)$ is derived. By inverting it, we obtain an asymptotic relation for the transition density of the process. We show that the error in this formula has the order $o(t^3)$ and, therefore, it gives a good approximation on small time intervals whose lengths depend on $\lambda$. Estimate of the accuracy of the approximation is analysed.  
\end{abstract}

\vskip 0.1cm


\vskip 0.2cm


\section{Introduction}

\numberwithin{equation}{section}

Among the great variety of the works devoted to random motions at finite speed in the Euclidean spaces $\Bbb R^m, \; m\ge 2,$ (see \cite{ghosh}, \cite{kol1,kol2,kol3,kol4,kol5,kol6}, \cite{mas}, \cite{sta1,sta2} for the Markovian case and \cite{lecaer1,lecaer2,lecaer3,letac}, \cite{pogor1,pogor2} for different non-Markovian cases), the Markov random flight in the three-dimensional Euclidean space $\Bbb R^3$ is, undoubtedly, the most difficult and hard to study. While in the low even-dimensional spaces $\Bbb R^2, \Bbb R^4$ and $\Bbb R^6$ one managed to obtain the distributions of the motions in an explicit form (see \cite{kol6}, \cite{kol5} and \cite{kol1}, respectively), in the important three-dimensional case only a few results are known.   

The absolutely continuous part of the transition density of the symmetric Markov random flight with unit speed in the Euclidean space $\Bbb R^3$ was presented in \cite[formulas (1.3) and (4.21) therein]{sta1}. It has an extremely complicated form of an integral with variable limits whose integrand involves inverse hyperbolic tangent function. This formula has so complicated form that cannot even be evaluated by means of standard computer environments. Moreover, the lack of the speed parameter in this formula impoverishes somewhat the model because it does not allow to study the limiting behaviour of the motion under various scaling conditions (under Kac's condition, for example). 
The presence of both parameters (i.e. the speed and the intensity of switchings) in any process of Markov random flight makes it, undoubtedly, the most adequate and realistic model for describing the finite-velocity diffusion in the Euclidean spaces. These parameters cannot be considered as independent because they are connected with each other through the time (namely, the speed is the distance passed {\it per unit of time} and the intensity is the mean number of switchings {\it per unit of time}). Another question concerning the density presented in \cite{sta1} is the infinite discontinuity at the origin $\bold 0\in\Bbb R^3$. While the infinite  discontinuity of the transition density on the border of the diffusion area is a quite natural property in some Euclidean spaces of low dimensions (see \cite{sta2,mas,kol6} for the Euclidean plane $\Bbb R^2$ and \cite[the second term of formulas (1.3) and (4.21)]{sta1}, \cite{kol4}, \cite[formula (3.12)]{kol2} in the space $\Bbb R^3$), the discontinuity at the origin looks somewhat strange and hard to explain. 

The difficulty of analysing the three-dimensional Markov random flight and, on the other hand, the great theoretical and applied importance of the problem of describing the finite-velocity diffusion in the space $\Bbb R^3$ suggest to look for other methods of studying this model. That is why various asymptotic theorems yielding a good approximation would be a fairly desirable aim of the research. Such asymptotic results could be obtained by using the characteristic functions technique. In the case of the three-dimensional symmetric Markov random flight some important results for its characteristic functions were obtained. In particular, the closed-form expression for the Laplace transform of the characteristic function was obtained by different methods in \cite[formulas (1.6) and (5.8)]{sta1} (for unit speed) and in \cite[formula (45)]{mas}, \cite[page 1054]{kol2} (for arbitrary speed). A general relation for the conditional characteristic functions of the three-dimensional symmetric Markov random flight conditioned by the number of changes of direction, was given in \cite[formula (3.8)]{kol2}. 

The key point in these formulas is the possibility of evaluating the inverse Laplace transforms of the powers of the inverse tangent functions in the complex right half-plane. This is the basic idea of deriving the series representations of the conditional characteristic functions corresponding to two and three changes of direction given in Section 3. Based on these representations, an asymptotic formula, as time $t\to 0$, for the unconditional characteristic function is obtained in Section 4 and the error in this formula has the order $o(t^3)$. The inverse Fourier transformation of the unconditional characteristic function yields an asymptotic formula for the transition density of the process which is presented in Section 5. This formula shows that the density is discontinuous on the border, but it is continuous at the origin $\bold 0\in\Bbb R^3$, as it must be. 
The unexpected and interesting peculiarity is that the conditional density corresponding to two changes of direction contains a term  having an infinite discontinuity on the border of the diffusion area. From this fact it follows that such conditional density is discontinuous itself on the border and this differs the 3D-model from its 2D-counterpart where only the conditional density of the single change of direction has an infinite discontinuity on the border. The error in the obtained asymptotic formula has the order $o(t^3)$. In Section 6 we estimate the accuracy of the asymptotic formula and show that it gives a good approximation on small time intervals whose lengths depend on the intensity of switchings. Finally, in Appendices we prove a series of auxiliary lemmas that have been used in our analysis.

\section{Description of the process and structure of distribution}

\numberwithin{equation}{section}

Consider the stochastic motion of a particle that, at the initial time instant $t=0$,  starts from the origin $\bold 0 = (0,0,0)$ of the Euclidean space $\Bbb R^3$ and moves with some constant speed $c$ (note that $c$ is treated as the constant norm of the velocity). The initial direction is a random three-dimensional vector with uniform distribution on the unit sphere 
$$S_1 = \left\{ \bold x=(x_1, x_2, x_3)\in \Bbb R^3: \; \Vert\bold x\Vert^2 = x_1^2+x_2+x_3^2=1 \right\} $$
The motion is controlled by a homogeneous Poisson process $N(t)$ of rate $\lambda>0$ as follows. At each Poissonian instant, the particle instantaneously takes on a new random direction distributed uniformly on $S_1$ independently of its previous motion and keeps moving with the same speed $c$ until the next Poisson event occurs, then it takes on a new random direction again and so on. 

Let $\bold X(t)=(X_1(t), X_2(t), X_3(t))$ be the particle's position at time $t>0$ which is referred to as the three-dimensional symmetric Markov random flight. At arbitrary time instant $t>0$ the particle, with probability 1, is located in the closed three-dimensional ball of radius $ct$ centred at the origin $\bold 0$:
$$\mathcal B_{ct} = \left\{ \bold x=(x_1, x_2, x_3)\in \Bbb R^3 : \; \Vert\bold x\Vert^2 = x_1^2+x_2^2+x_3^2\le c^2t^2 \right\} .$$

Consider the probability distribution function 
$$\Phi(\bold x, t) = \text{Pr} \left\{ \bold X(t)\in d\bold x \right\}, \qquad \bold x\in\mathcal B_{ct}, \quad t\ge 0,$$ 
of the process $\bold X(t)$, where $d\bold x$ is the infinitesimal element in the space $\Bbb R^3$. For arbitrary fixed $t>0$, the distribution $\Phi(\bold x, t)$ consists of two components. 

The singular component corresponds to the case when no Poisson events occur on the time interval $(0,t)$ and it is concentrated on the sphere
$$S_{ct} =\partial\mathcal B_{ct} = \left\{ \bold x=(x_1, x_2, x_3)\in \Bbb R^3: \; \Vert\bold x\Vert^2 = x_1^2+x_2^2+x_3^2=c^2t^2 \right\} .$$
In this case, at time instant $t$, the particle is located on the sphere $S_{ct}$ and the probability of this event is
$$\text{Pr} \left\{ \bold X(t)\in S_{ct} \right\} = e^{-\lambda t} .$$

If at least one Poisson event occurs on the time interval $(0, t)$, then the particle is located strictly inside the ball $\mathcal B_{ct}$ and the probability of this event is
\begin{equation}\label{struc1}
\text{Pr} \left\{ \bold X(t)\in \text{int} \; \mathcal B_{ct} \right\} = 1 - e^{-\lambda t} .
\end{equation}
The part of the distribution $\Phi(\bold x, t)$ corresponding to this case is concentrated in the interior 
$$\text{int} \; \mathcal B_{ct} = \left\{ \bold x=(x_1, x_2, x_3)\in \Bbb R^3: \; \Vert\bold x\Vert^2 = x_1^2+x_2^2+x_3^2<c^2t^2 \right\}$$
of the ball $\mathcal B_{ct}$ and forms its absolutely continuous component. 

Let $p(\bold x,t), \; \bold x\in\mathcal B_{ct} , \; t>0,$ be the density of distribution $\Phi(\bold x,t)$. It has the form 
\begin{equation}\label{struc2}
p(\bold x,t) = p^{(s)}(\bold x,t) + p^{(ac)}(\bold x,t) , \qquad \bold x\in\mathcal B_{ct}, \quad t>0,
\end{equation}
where $p^{(s)}(\bold x,t)$ is the density (in the sense of generalized functions) of the singular component of $\Phi(\bold x, t)$ concentrated on the sphere $S_{ct}$ and $p^{(ac)}(\bold x,t)$ is the density of the absolutely continuous component of $\Phi(\bold x, t)$ concentrated in $\text{int} \; \mathcal B_{ct}$.

The singular part of density (\ref{struc2}) is given by the formula: 
\begin{equation}\label{densS}
p^{(s)}(\bold x,t) =  \frac{e^{-\lambda t}}{4\pi c^2t^2} \; \delta(c^2t^2-\Vert\bold x\Vert^2) , \qquad t>0,
\end{equation}
where $\delta(x)$ is the Dirac delta-function. 

The absolutely continuous part of density (\ref{struc2}) has the form:
\begin{equation}\label{densAC}
p^{(ac)}(\bold x,t) = f^{(ac)}(\bold x,t) \Theta(ct-\Vert\bold x\Vert) , \qquad t>0,
\end{equation}
where $f^{(ac)}(\bold x,t)$ is some positive function absolutely continuous in $\text{int} \; \mathcal B_{ct}$ and $\Theta(x)$ is the Heaviside unit-step function given by 
\begin{equation}\label{heaviside}
\Theta(x) = \left\{ \aligned 1, \qquad  & \text{if} \; x>0,\\
                               0, \qquad & \text{if} \; x\le 0.
\endaligned \right.
\end{equation}

Asymptotic behaviour of the transition density (\ref{struc2}) on small time intervals is the main subject of this research. Since its singular part is explicitly given by (\ref{densS}), then our efforts are mostly concentrated on deriving the respective asymptotic formulas for the absolutely continuous component (\ref{densAC}) of the density. Our main tool is the characteristic functions technique because, as it was mentioned above, some closed-form expressions for the characteristic functions (both conditional and unconditional ones) of the three-dimensional symmetric Markov random flight $\bold X(t)$ are known.

\section{Conditional characteristic functions}

\numberwithin{equation}{section}

In this section we obtain the series representations of the conditional characteristic functions corresponding to two and three changes of direction. These formulas are the basis for our further analysis leading to asymptotic relations for the unconditional characteristic function and the transition density of the three-dimensional symmetric Markov random flight $\bold X(t)$ on small time intervals. The main result of this section is given by the following theorem. 

\bigskip 

{\bf Theorem 1.} {\it The conditional characteristic functions $H_2(\boldsymbol\alpha,t)$ and $H_3(\boldsymbol\alpha,t)$ 
corresponding to two and three changes of direction are given, respectively, by the formulas:}
\begin{equation}\label{char2}
\aligned 
H_2(\boldsymbol\alpha,t) & = \sum_{k=0}^{\infty} \frac{(ct\Vert\boldsymbol\alpha\Vert)^{k-1}}{2^{k-1} \; k! \; (2k+1)^2} \\ 
& \qquad \times \;  _5F_4\left( 1,1,1,-k,-k-\frac{1}{2}; \; -k+\frac{1}{2}, -k+\frac{1}{2}, \frac{3}{2}, 2; \; 1 \right) J_{k+1}(ct\Vert\boldsymbol\alpha\Vert) , 
\endaligned
\end{equation}
\begin{equation}\label{char3}
\aligned 
H_3(\boldsymbol\alpha,t) & = 3\pi^{3/2} \; \sum_{k=0}^{\infty} \frac{\gamma_k \; (ct\Vert\boldsymbol\alpha\Vert)^{k-3/2}}{2^{k+3/2} \; (k+1)!} \;  J_{k+3/2}(ct\Vert\boldsymbol\alpha\Vert) ,
\endaligned
\end{equation}
$$\boldsymbol\alpha = (\alpha_1,\alpha_2,\alpha_3) \in\Bbb R^3, \qquad \Vert\boldsymbol\alpha\Vert = \sqrt{\alpha_1^2+\alpha_2^2+\alpha_3^2}, \qquad t>0, $$ 
{\it where $J_{\nu}(z)$ is Bessel function, $_5F_4(a_1,a_2,a_3,a_4,a_5; \; b_1,b_2,b_3,b_4; \; z)$ 
is the generalized hypergeometric function given by} (\ref{hypergeom54}) {\it (see below) and the coefficients $\gamma_k$ are given by the formula} 
\begin{equation}\label{coef1}
\gamma_k = \frac{1}{k+2} \sum_{l=0}^k \frac{l! \; (k-l)!}{(l+1) \; \Gamma\left( l+\frac{3}{2} \right) \; \Gamma\left( k-l+\frac{3}{2} \right)} , \qquad k\ge 0. 
\end{equation}

\vskip 0.2cm

{\it Proof.} 
It was proved in \cite[formula (3.8)]{kol2} that, for arbitrary $t>0$, the characteristic function $H_n(\boldsymbol\alpha,t)$ (that is, Fourier transform 
$\mathcal F_{\bold x}$ with respect to spatial variable $\bold x$)  of the conditional density $p_n(\bold x,t)$ of the three-dimensional Markov random flight 
$\bold X(t)$ corresponding to $n$ changes of directions is given by the formula
\begin{equation}\label{eq1}
H_n(\boldsymbol\alpha,t) = \mathcal F_{\bold x}[p_n(\bold x,t)](\boldsymbol\alpha) = \frac{n!}{t^n} (c\Vert\boldsymbol\alpha\Vert)^{-(n+1)} \mathcal L_s^{-1} \left[ \left( \text{arctg} \frac{c\Vert\boldsymbol\alpha\Vert}{s} \right)^{n+1} \right](t) , 
\end{equation}
$$n\ge 1, \qquad \boldsymbol\alpha\in\Bbb R^3, \quad s\in\Bbb C^+ ,$$
where $\mathcal L_s^{-1}$ is the inverse Laplace transformation with respect to complex variable $s$ and $\Bbb C^+ = \{ s\in\Bbb C : \text{Re} \; s >0 \}$ is the right half-plane of the complex plane $\Bbb C$. In particular, in the case of two changes of directions $n=2$, formula (\ref{eq1}) yields: 
\begin{equation}\label{eq2}
H_2(\boldsymbol\alpha,t) = \mathcal F_{\bold x}[p_2(\bold x,t)](\boldsymbol\alpha) = \frac{2!}{t^2} (c\Vert\boldsymbol\alpha\Vert)^{-3} \mathcal L_s^{-1} \left[ \left( \text{arctg} \frac{c\Vert\boldsymbol\alpha\Vert}{s} \right)^3 \right](t) , \qquad \boldsymbol\alpha\in\Bbb R^3, \quad s\in\Bbb C^+ . 
\end{equation}
Applying Lemma B3 of the Appendix B to the power of inverse tangent function in (\ref{eq2}), we obtain: 
\begin{equation}\label{eq3}
\aligned 
& H_2(\boldsymbol\alpha,t) \\ 
& = \frac{2}{t^2} (c\Vert\boldsymbol\alpha\Vert)^{-3} \mathcal L_s^{-1} \biggl[ \frac{1}{\sqrt{\pi}} \left( \frac{c\Vert\boldsymbol\alpha\Vert}{\sqrt{s^2+(c\Vert\boldsymbol\alpha\Vert)^2}} \right)^3  \; \sum_{k=0}^{\infty} \frac{\Gamma\left( k+\frac{1}{2} \right)}{k! \; (2k+1)} \\ 
& \qquad \times  _5F_4\left( 1,1,1,-k,-k-\frac{1}{2}; \; -k+\frac{1}{2}, -k+\frac{1}{2}, \frac{3}{2}, 2; \; 1 \right) \left( \frac{(c\Vert\boldsymbol\alpha\Vert)^2}{s^2 + (c\Vert\boldsymbol\alpha\Vert)^2} \right)^k \biggr](t) \\
& = \frac{2}{\sqrt{\pi} \; t^2} \; \sum_{k=0}^{\infty} \frac{\Gamma\left( k+\frac{1}{2} \right)}{k! \; (2k+1)} \; (c\Vert\boldsymbol\alpha\Vert)^{2k} \\ 
& \qquad \times \;  _5F_4\left( 1,1,1,-k,-k-\frac{1}{2}; \; -k+\frac{1}{2}, -k+\frac{1}{2}, \frac{3}{2}, 2; \; 1 \right) \mathcal L_s^{-1} \biggl[ \frac{1}{\left( s^2 + (c\Vert\boldsymbol\alpha\Vert)^2 \right)^{k+3/2}} \biggr](t) .
\endaligned
\end{equation}
Note that evaluating the inverse Laplace transformation of each term of the series separately is justified because it converges uniformly in $s$ everywhere in $\Bbb C^+$ and the complex functions 
$\left( s^2 + (c\Vert\boldsymbol\alpha\Vert)^2 \right)^{-(k+3/2)}, \; k\ge 0,$ are holomorphic and do not have any singular points in this half-plane. Moreover, each of these functions contains the inversion complex variable $s\in\Bbb C^+$ in a negative power and behaves like $s^{-(2k+3)}$, as $|s|\to +\infty$, and, therefore, all these complex functions rapidly tend to zero at infinity. 

According to \cite[Table 8.4-1, formula 57]{korn}, we have 
$$\mathcal L_s^{-1} \biggl[ \frac{1}{\left( s^2 + (c\Vert\boldsymbol\alpha\Vert)^2 \right)^{k+3/2}} \biggr](t) = \frac{\sqrt{\pi}}{\Gamma\left( k+\frac{3}{2} \right)} 
\left( \frac{t}{2c\Vert\boldsymbol\alpha\Vert} \right)^{k+1} J_{k+1}(ct\Vert\boldsymbol\alpha\Vert) .$$
Substituting this into (\ref{eq3}), after some simple calculations we obtain (\ref{char2}).

For $n=3$, formula (\ref{eq1}) yields: 
\begin{equation}\label{eq4}
H_3(\boldsymbol\alpha,t) = \mathcal F_{\bold x}[p_3(\bold x,t)](\boldsymbol\alpha) = \frac{3!}{t^3} (c\Vert\boldsymbol\alpha\Vert)^{-4} \mathcal L_s^{-1} \left[ \left( \text{arctg} \frac{c\Vert\boldsymbol\alpha\Vert}{s} \right)^4 \right](t) , \qquad \boldsymbol\alpha\in\Bbb R^3, \quad s\in\Bbb C^+ .
\end{equation}
Applying Lemma B4 of the Appendix B to the power of inverse tangent function in (\ref{eq4}) and taking into account that
$$\mathcal L_s^{-1} \biggl[ \frac{1}{\left( s^2 + (c\Vert\boldsymbol\alpha\Vert)^2 \right)^{k+2}} \biggr](t) = 
\frac{\sqrt{\pi}}{(k+1)!} 
\left( \frac{t}{2c\Vert\boldsymbol\alpha\Vert} \right)^{k+3/2} J_{k+3/2}(ct\Vert\boldsymbol\alpha\Vert) ,$$
we obtain: 
\begin{equation}\label{eq5}
\aligned 
H_3(\boldsymbol\alpha,t) & = \frac{3\pi}{t^3} \sum_{k=0}^{\infty} \gamma_k \; (c\Vert\boldsymbol\alpha\Vert)^{2k} \; \mathcal L_s^{-1} \biggl[ \frac{1}{\left( s^2 + (c\Vert\boldsymbol\alpha\Vert)^2 \right)^{k+2}} \biggr](t) \\
& = 3\pi^{3/2} \; \sum_{k=0}^{\infty} \frac{\gamma_k \; (ct\Vert\boldsymbol\alpha\Vert)^{k-3/2}}{2^{k+3/2} \; (k+1)!} \;  J_{k+3/2}(ct\Vert\boldsymbol\alpha\Vert) ,
\endaligned
\end{equation}
where the coefficients $\gamma_k$ are given by (\ref{coef1}). The theorem is proved. $\square$ 

\bigskip 

{\it Remark 1.} The series in formulas (\ref{char2}) and (\ref{char3}) are convergent for any fixed $t>0$, however this convergence is not uniform in $\Vert\boldsymbol\alpha\Vert$. Therefore, we cannot invert each term of these series separately. Moreover, one can see  that the inverse Fourier transform of each term does not exist for $k\ge 2$. Thus, while there exist the inverse Fourier transforms of the whole series (\ref{char2}) and (\ref{char3}), it is impossible to invert their terms separately and, therefore, we cannot obtain closed-form expressions for the respective conditional densities. These formulas can, nevertheless, be used for obtaining the important  asymptotic relations and this is the main subject of the next sections.

\section{Asymptotic formula for characteristic function}

\numberwithin{equation}{section}

Using the results of the previous section, we can now present an asymptotic relation on small time intervals for the characteristic function
$$H(\boldsymbol\alpha,t) = e^{-\lambda t} \sum_{k=0}^{\infty} \frac{(\lambda t)^k}{k!} H_k(\boldsymbol\alpha,t) $$
of the three-dimensional symmetric Markov random flight, where $H_k(\boldsymbol\alpha,t), \; k\ge 0,$ are the conditional characteristic functions corresponding to $k$ changes of direction. This result is given by the following theorem. 

\bigskip 

{\bf Theorem 2.} {\it For the characterictic function $H(\boldsymbol\alpha,t), \; t>0,$ of the three-dimensional Markov random flight $\bold X(t)$ the following asymptotic formula holds:} 
\begin{equation}\label{eq6}
\aligned 
H(\boldsymbol\alpha,t) & = e^{-\lambda t} \biggl\{ \frac{\sin{(ct\Vert\boldsymbol\alpha\Vert)}}{ct\Vert\boldsymbol\alpha\Vert} + \frac{\lambda}{c^2 t \Vert\boldsymbol\alpha\Vert^2} \biggl[ \sin{(ct\Vert\boldsymbol\alpha\Vert)} \text{Si}(2ct\Vert\boldsymbol\alpha\Vert) + \cos{(ct\Vert\boldsymbol\alpha\Vert)} \text{Ci}(2ct\Vert\boldsymbol\alpha\Vert) \biggr] \\ 
& \qquad\qquad\qquad + \frac{\lambda^2 t}{c\Vert\boldsymbol\alpha\Vert} J_1(ct\Vert\boldsymbol\alpha\Vert) + \frac{\lambda^3 \; \sqrt{\pi} \; t^{3/2}}{(2 c\Vert\boldsymbol\alpha\Vert)^{3/2}} \; J_{3/2}(ct\Vert\boldsymbol\alpha\Vert)\biggr\} + o(t^3) , 
\endaligned
\end{equation}
$$\boldsymbol\alpha = (\alpha_1,\alpha_2,\alpha_3) \in\Bbb R^3, \qquad \Vert\boldsymbol\alpha\Vert = \sqrt{\alpha_1^2+\alpha_2^2+\alpha_3^2}, \qquad t>0, $$ 
{\it where} $\text{Si}(z)$ {\it and} $\text{Ci}(z)$ {\it are the incomplete integral sine and cosine, respectively, given by the formulas:}
$$\text{Si}(x) = \int_0^x \frac{\sin\xi}{\xi} \; d\xi, \qquad
\text{Ci}(x) = \int_0^x \frac{\cos\xi -1}{\xi} \; d\xi.$$

\vskip 0.2cm 

{\it Proof.} We have:  
$$H(\boldsymbol\alpha,t) = e^{-\lambda t} \biggl[ H_0(\boldsymbol\alpha,t) + \lambda t H_1(\boldsymbol\alpha,t) + \frac{(\lambda t)^2}{2!} H_2(\boldsymbol\alpha,t) + \frac{(\lambda t)^3}{3!} H_3(\boldsymbol\alpha,t) + \sum_{k=4}^{\infty} \frac{(\lambda t)^k}{k!} H_k(\boldsymbol\alpha,t) \biggr] .$$ 
Since all the conditional characteristic functions are uniformly bounded in both variables, that is, $\vert H_k(\boldsymbol\alpha,t) \vert \le 1, \; \boldsymbol\alpha\in\Bbb R, \; t\ge 0, \; k\ge 0,$ then
$$\sum_{k=4}^{\infty} \frac{(\lambda t)^k}{k!} H_k(\boldsymbol\alpha,t) = o(t^3)$$
and, therefore, 
\begin{equation}\label{eq7}
H(\boldsymbol\alpha,t) = e^{-\lambda t} \biggl[ H_0(\boldsymbol\alpha,t) + \lambda t H_1(\boldsymbol\alpha,t) + \frac{(\lambda t)^2}{2!} H_2(\boldsymbol\alpha,t) + \frac{(\lambda t)^3}{3!} H_3(\boldsymbol\alpha,t) + o(t^3) \biggr].  
\end{equation}
In view of (\ref{char2}), we have:  
$$\aligned 
& \frac{(\lambda t)^2}{2!} H_2(\boldsymbol\alpha,t) \\
& = \lambda^2 \biggl[ \frac{t}{c\Vert\boldsymbol\alpha\Vert} J_1(ct\Vert\boldsymbol\alpha\Vert) \\
& \quad + \sum_{k=1}^{\infty} \frac{(c\Vert\boldsymbol\alpha\Vert)^{k-1} t^{k+1}}{2^k \; k! \; (2k+1)^2}  
\;  _5F_4\left( 1,1,1,-k,-k-\frac{1}{2}; \; -k+\frac{1}{2}, -k+\frac{1}{2}, \frac{3}{2}, 2; \; 1 \right) J_{k+1}(ct\Vert\boldsymbol\alpha\Vert) \biggr] . 
\endaligned$$
From the asymptotic formula 
\begin{equation}\label{asbes}
J_{\nu}(z) = \frac{z^{\nu}}{2^{\nu} \; \Gamma(\nu+1)} + o(z^{\nu+1}) , \qquad \nu\ge 0,  
\end{equation}
we get 
$$J_{k+1}(ct\Vert\boldsymbol\alpha\Vert) = \frac{(ct\Vert\boldsymbol\alpha\Vert)^{k+1}}{2^{k+1} (k+1)!} + o(t^{k+2})$$
and, therefore,
$$\sum_{k=1}^{\infty} \frac{(c\Vert\boldsymbol\alpha\Vert)^{k-1} t^{k+1}}{2^k \; k! \; (2k+1)^2}  
\;  _5F_4\left( 1,1,1,-k,-k-\frac{1}{2}; \; -k+\frac{1}{2}, -k+\frac{1}{2}, \frac{3}{2}, 2; \; 1 \right) J_{k+1}(ct\Vert\boldsymbol\alpha\Vert) = o(t^3) .$$
Thus, we obtain the following asymptotic relation: 
\begin{equation}\label{eq8}
\frac{(\lambda t)^2}{2!} H_2(\boldsymbol\alpha,t) = \frac{\lambda^2 t}{c\Vert\boldsymbol\alpha\Vert} J_1(ct\Vert\boldsymbol\alpha\Vert) + o(t^3) .
\end{equation}

Similarly, according to (\ref{char3}), we have:
$$\aligned
\frac{(\lambda t)^3}{3!} H_3(\boldsymbol\alpha,t) & = \lambda^3 \pi^{3/2} \biggl[ \frac{\gamma_0 (c\Vert\boldsymbol\alpha\Vert)^{-3/2} t^{3/2}}{2^{5/2}} J_{3/2}(ct\Vert\boldsymbol\alpha\Vert) \\
& \qquad\qquad\qquad + \sum_{k=1}^{\infty} \frac{\gamma_k \; (c\Vert\boldsymbol\alpha\Vert)^{k-3/2} t^{k+3/2}}{2^{k+5/2} \; (k+1)!} \;  J_{k+3/2}(ct\Vert\boldsymbol\alpha\Vert) \biggr] .
\endaligned$$
In view of (\ref{asbes}), we have
$$J_{k+3/2}(ct\Vert\boldsymbol\alpha\Vert) = \frac{(ct\Vert\boldsymbol\alpha\Vert)^{k+3/2}}{2^{k+3/2} \; \Gamma\left( k+\frac{5}{2}  \right)} + o(t^{k+5/2})$$
and, therefore, 
$$\sum_{k=1}^{\infty} \frac{\gamma_k \; (c\Vert\boldsymbol\alpha\Vert)^{k-3/2} t^{k+3/2}}{2^{k+5/2} \; (k+1)!} \;  J_{k+3/2}(ct\Vert\boldsymbol\alpha\Vert) = o(t^4) .$$
Thus, taking into account that $\gamma_0 = 2/\pi$ (see (\ref{coef1})), we arrive at the formula: 
\begin{equation}\label{eq9}
\frac{(\lambda t)^3}{3!} H_3(\boldsymbol\alpha,t) = \frac{\lambda^3 \; \sqrt{\pi} \; t^{3/2}}{(2 c\Vert\boldsymbol\alpha\Vert)^{3/2}} 
\; J_{3/2}(ct\Vert\boldsymbol\alpha\Vert) + o(t^4) .
\end{equation}

Since (see \cite[formula (3.11)]{kol2})  
$$\lambda t H_1(\boldsymbol\alpha,t) = \frac{\lambda}{c^2 t \Vert\boldsymbol\alpha\Vert^2} \biggl[ \sin{(ct\Vert\boldsymbol\alpha\Vert)} \text{Si}(2ct\Vert\boldsymbol\alpha\Vert) + \cos{(ct\Vert\boldsymbol\alpha\Vert)} \text{Ci}(2ct\Vert\boldsymbol\alpha\Vert) \biggr]$$
and 
$$H_0(\boldsymbol\alpha,t) = \frac{\sin{(ct\Vert\boldsymbol\alpha\Vert)}}{ct\Vert\boldsymbol\alpha\Vert}$$
(that is, characteristic function of the uniform distribution on the surface of the three-dimensional sphere of radius $ct$), then by substituting these formulas, as well as (\ref{eq8}) and (\ref{eq9}) into (\ref{eq7}), we finally obtain asymptotic relation (\ref{eq6}). The theorem is completely proved. $\square$

\section{Asymptotic relation for the transition density}

\numberwithin{equation}{section}

Asymptotic formula (\ref{eq6}) for the unconditional characteristic function enables us to obtain the respective asymptotic relation for the transition density of the process $\bold X(t)$. This result is given by the following theorem. 

\bigskip

{\bf Theorem 3.} {\it For the transition density $p(\bold x,t), \; t>0,$ of the three-dimensional Markov random flight $\bold X(t)$ the following asymptotic relation holds:} 
\begin{equation}\label{dens1}
\aligned 
p(\bold x,t) & = \frac{e^{-\lambda t}}{4\pi (ct)^2} \; \delta(c^2t^2-\Vert\bold x\Vert^2) + e^{-\lambda t} \biggl[ \frac{\lambda}{4\pi c^2 t \Vert\bold x\Vert} \ln\left( \frac{ct+\Vert\bold x\Vert}{ct-\Vert\bold x\Vert} \right) \\ 
& \qquad\qquad  + \frac{\lambda^2}{2\pi^2 c^2 \; \sqrt{c^2t^2-\Vert\bold x\Vert^2}} + \frac{\lambda^3}{8\pi c^3} \biggr] \Theta(ct-\Vert\bold x\Vert) + o(t^3) , 
\endaligned
\end{equation}
$$\bold x = (x_1, x_2, x_3) \in\Bbb R^3, \qquad \Vert\bold x\Vert = \sqrt{x_1^2+x_2^2+x_3^2}, \qquad t>0.$$

\vskip 0.2cm 

{\it Proof.} Applying the inverse Fourier transformation $\mathcal F_{\boldsymbol\alpha}^{-1}$ to both sides of (\ref{eq6}), we have: 
\begin{equation}\label{dens2}
\aligned 
p(\bold x,t) & = e^{-\lambda t} \biggl\{ \mathcal F_{\boldsymbol\alpha}^{-1} \biggl[  \frac{\sin{(ct\Vert\boldsymbol\alpha\Vert)}}{ct\Vert\boldsymbol\alpha\Vert} \biggr](\bold x) \\
& \qquad\quad + \mathcal F_{\boldsymbol\alpha}^{-1} \biggl[ \frac{\lambda}{c^2 t \Vert\boldsymbol\alpha\Vert^2} \biggl( \sin{(ct\Vert\boldsymbol\alpha\Vert)} \text{Si}(2ct\Vert\boldsymbol\alpha\Vert) + \cos{(ct\Vert\boldsymbol\alpha\Vert)} \text{Ci}(2ct\Vert\boldsymbol\alpha\Vert) \biggr) \biggr](\bold x) \\ 
& \qquad\quad + \mathcal F_{\boldsymbol\alpha}^{-1} \biggl[ \frac{\lambda^2 t}{c\Vert\boldsymbol\alpha\Vert} J_1(ct\Vert\boldsymbol\alpha\Vert) \biggr](\bold x) \\
& \qquad\quad + \mathcal F_{\boldsymbol\alpha}^{-1} \biggl[ \frac{\lambda^3 t}{2(c\Vert\boldsymbol\alpha\Vert)^2} \left( \frac{\sin{(ct\Vert\boldsymbol\alpha\Vert)}}{ct\Vert\boldsymbol\alpha\Vert} - \cos{(ct\Vert\boldsymbol\alpha\Vert)} \right) \biggr](\bold x) \biggr\} + o(t^3) . 
\endaligned
\end{equation}
Note that here we have used the fact that, due to the continuity of the inverse Fourier transformation, the asymptotic formula $\mathcal F_{\boldsymbol\alpha}^{-1} \bigl[ o(t^3) \bigr](\bold x) = o(t^3)$ holds. 

Let us evaluate separately the inverse Fourier transforms on the right-hand side of (\ref{dens2}). The first one is well known (see \cite[page 1051]{kol2}): 
\begin{equation}\label{dens3}
\mathcal F_{\boldsymbol\alpha}^{-1} \biggl[  \frac{\sin{(ct\Vert\boldsymbol\alpha\Vert)}}{ct\Vert\boldsymbol\alpha\Vert} \biggr](\bold x) = \frac{1}{4\pi (ct)^2} \; \delta(c^2t^2-\Vert\bold x\Vert^2) 
\end{equation}
that is the uniform density concentrated on the surface of the sphere $S_{ct} \in\Bbb R^3$ of radius $ct$ centred at the origin $\bold 0\in\Bbb R^3$.  

The second Fourier transform on the right-hand side of (\ref{dens2}) is also well known (see \cite[page 64, the Theorem]{kol4} or 
\cite[formulas (3.11) and (3.12)]{kol2}): 
\begin{equation}\label{dens4}
\aligned 
\mathcal F_{\boldsymbol\alpha}^{-1} & \biggl[ \frac{\lambda}{c^2 t \Vert\boldsymbol\alpha\Vert^2} \biggl( \sin{(ct\Vert\boldsymbol\alpha\Vert)} \text{Si}(2ct\Vert\boldsymbol\alpha\Vert) + \cos{(ct\Vert\boldsymbol\alpha\Vert)} \text{Ci}(2ct\Vert\boldsymbol\alpha\Vert) \biggr) \biggr](\bold x) \\
& \hskip 4cm = \frac{\lambda}{4\pi c^2 t \Vert\bold x\Vert} \ln\left( \frac{ct+\Vert\bold x\Vert}{ct-\Vert\bold x\Vert} \right) \; \Theta(ct-\Vert\bold x\Vert) .
\endaligned  
\end{equation}

Applying the Hankel inversion formula, we have for the third Fourier transform on the right-hand side of (\ref{dens2}):
$$\frac{\lambda^2 t}{c} \; \mathcal F_{\boldsymbol\alpha}^{-1} \biggl[ \Vert\boldsymbol\alpha\Vert^{-1} \;  J_1(ct\Vert\boldsymbol\alpha\Vert) \biggr](\bold x) = \frac{\lambda^2 t}{c} \; (2\pi)^{-3/2} \Vert\bold x\Vert^{-1/2} \int_0^{\infty} J_{1/2}(\Vert\bold x\Vert \xi) \;  \xi^{3/2} \; \xi^{-1} J_1(ct\xi) \; d\xi .$$
Taking into account that 
\begin{equation}\label{bessin}
J_{1/2}(z) = \sqrt{\frac{2}{\pi z}} \; \sin{z},
\end{equation}
and applying \cite[formula 2.12.15(2)]{pbm}, we have: 
\begin{equation}\label{dens5}
\aligned 
\frac{\lambda^2 t}{c} \; \mathcal F_{\boldsymbol\alpha}^{-1} \biggl[ \Vert\boldsymbol\alpha\Vert^{-1} \;  J_1(ct\Vert\boldsymbol\alpha\Vert) \biggr](\bold x) & = \frac{\lambda^2 t}{2\pi^2 c \Vert\bold x\Vert} \int_0^{\infty} \sin{(\Vert\bold x\Vert \xi)} \; J_1(ct\xi) 
\; d\xi \\ 
& = \frac{\lambda^2 t}{2\pi^2 c \Vert\bold x\Vert} \; (c^2t^2-\Vert\bold x\Vert^2)^{-1/2} \; 
\left( \frac{\Vert\bold x\Vert}{ct} \right) \; \Theta(ct-\Vert\bold x\Vert)\\ 
& = \frac{\lambda^2}{2\pi^2 c^2 \; \sqrt{c^2t^2-\Vert\bold x\Vert^2}} \; \Theta(ct-\Vert\bold x\Vert) . 
\endaligned  
\end{equation}
This is a fairly unexpected result showing that the conditional density $p_2(\bold x,t)$ corresponding to two changes of direction has an infinite discontinuity on the border of the three-dimensional ball $\mathcal B_{ct}$. This property is similar to that of the conditional density $p_1(\bold x,t)$ corresponding to the single change of direction (for the respective joint density see  (\ref{dens4})).

Applying the Hankel inversion formula and taking into account (\ref{bessin}), we have for the fourth term 
on the right-hand side of (\ref{dens2}): 
$$\aligned 
\frac{\lambda^3 \; \sqrt{\pi} \; t^{3/2}}{(2 c)^{3/2}} \; & \mathcal F_{\boldsymbol\alpha}^{-1} \biggl[ \Vert\boldsymbol\alpha\Vert^{-3/2} \; J_{3/2}(ct\Vert\boldsymbol\alpha\Vert) \biggr](\bold x) \\
& = \frac{\lambda^3 \; \sqrt{\pi} \; t^{3/2}}{(2 c)^{3/2}} \; (2\pi)^{-3/2} \Vert\bold x\Vert^{-1/2} \int_0^{\infty} J_{1/2}(\Vert\bold x\Vert \xi) \;  \xi^{3/2} \; \xi^{-3/2} J_{3/2}(ct\xi) \; d\xi \\ 
& = \frac{\lambda^3 \; \sqrt{2} \; t^{3/2}}{8c^{3/2} \; \pi\sqrt{\pi} \; \Vert\bold x\Vert} \int_0^{\infty} \xi^{-1/2} \; \sin{(\Vert\bold x\Vert \xi)} \; J_{3/2}(ct\xi) \; d\xi . 
\endaligned$$
Using \cite[formula 6.699(1)]{gr}, we obtain: 
\begin{equation}\label{dens6}
\aligned 
\frac{\lambda^3 \; \sqrt{\pi} \; t^{3/2}}{(2 c)^{3/2}} \; & \mathcal F_{\boldsymbol\alpha}^{-1} \biggl[ \Vert\boldsymbol\alpha\Vert^{-3/2} \; J_{3/2}(ct\Vert\boldsymbol\alpha\Vert) \biggr](\bold x) \\ 
& = \frac{\lambda^3 \; \sqrt{2} \; t^{3/2}}{8c^{3/2} \; \pi\sqrt{\pi} \; \Vert\bold x\Vert} \; \frac{2^{-1/2} \; \sqrt{\pi} \; \Vert\bold x\Vert \; (ct)^{-3/2}}{\Gamma(1)} \; \Theta(ct-\Vert\bold x\Vert) \\ 
& = \frac{\lambda^3}{8\pi c^3} \; \Theta(ct-\Vert\bold x\Vert) .
\endaligned  
\end{equation}

Substituting now (\ref{dens3}), (\ref{dens4}), (\ref{dens5}) and (\ref{dens6}) into (\ref{dens2}) we arrive at (\ref{dens1}).
The theorem is proved. $\square$ 

\begin{center}
\begin{figure}[htbp]
\centerline{\includegraphics[width=10cm,height=8cm]{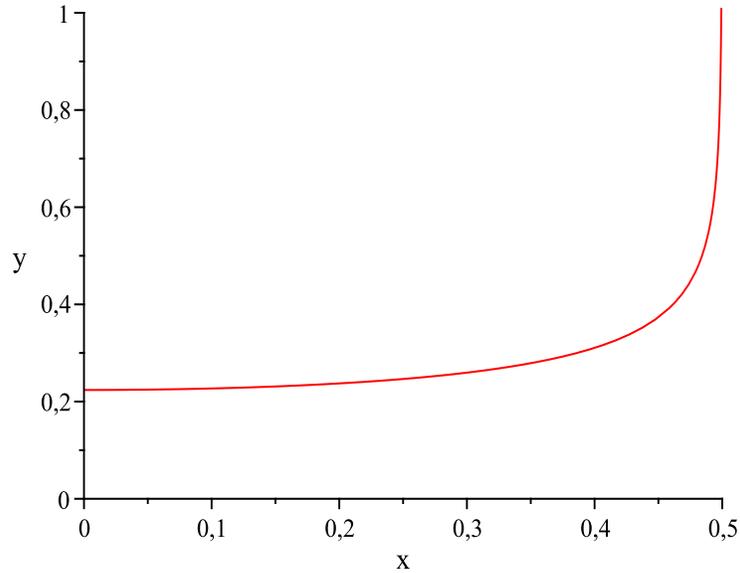}}
\caption{\it The shape of the absolutely continuous part of density (\ref{dens1}) at instant $t=0.1$ $\quad$   
(for $c=5, \; \lambda =2$) on the interval $\Vert\bold x\Vert \in [0, \; 0.5)$} 
\end{figure}
\end{center}
\vskip-1cm

The shape of the absolutely continuous part of density (\ref{dens1}) at time instant $t=0.1$ (for $c=5, \; \lambda=2$) on the interval $\Vert\bold x\Vert \in [0, \; 0.5)$ is plotted in Fig. 1. The error in these calculations does not exceed 0.001.

We see that the density increases slowly as the distance $\Vert\bold x\Vert$ from the origin $\bold 0\in\Bbb R^3$ grows, while near the border this growth becomes explosive. From this fact it follows that, for small time $t$, the greater part of the density is concentrated outside the neighbourhood of the origin $\bold 0\in\Bbb R^3$ and this feature of the three-dimensional Markov random flight is quite similar to that of  its two-dimensional counterpart. The infinite discontinuity of the density on the border $\Vert\bold x\Vert=ct$ is also similar to the analogous property of the two-dimensional Markov random flight (see, for comparison, \cite[formula (20) and Figure 2 therein]{kol6}). Note  that density (\ref{dens1}) is continuous at the origin, as it must be.  

\bigskip 

{\it Remark 2.} Using (\ref{dens1}), we can derive an asymptotic formula, as $t\to 0$, for the probability of being in a subball $\mathcal B_r$ of some radius $r<ct$ centred at the origin $\bold 0\in\Bbb R^3$. Applying \cite[formula 4.642]{gr} and \cite[formula 1.513(1)]{gr}, we have: 
\begin{equation}\label{dens7}
\aligned 
\int\limits_{\mathcal B_r} \frac{1}{\Vert\bold x\Vert} \; \ln\left( \frac{ct+\Vert\bold x\Vert}{ct-\Vert\bold x\Vert} \right) \; d\bold x & = 
\frac{2\pi^{3/2}}{\Gamma\left( \frac{3}{2} \right)} \int_0^r \xi^2 \; \frac{1}{\xi} \; \ln\left( \frac{ct+\xi}{ct-\xi} \right) \; d\xi \\ 
& = 4\pi (ct)^2 \int_0^{r/(ct)} z \; \ln\left( \frac{1+z}{1-z} \right) \; dz \\ 
& = 8\pi (ct)^2 \sum_{k=1}^{\infty} \frac{1}{2k-1} \int_0^{r/(ct)} z^{2k} \; dz \\ 
& = 8\pi r ct \sum_{k=1}^{\infty} \frac{1}{4k^2-1} \; \left( \frac{r^2}{c^2t^2} \right)^k .
\endaligned  
\end{equation}
This series can be expressed through the special Lerch $\psi$-function. 

Applying again \cite[formula 4.642]{gr}, we get: 
\begin{equation}\label{dens8}
\aligned 
\int\limits_{\mathcal B_r} \frac{d\bold x}{\sqrt{c^2t^2-\Vert\bold x\Vert}} & = \frac{2\pi^{3/2}}{\Gamma\left( \frac{3}{2} \right)} \int_0^r \frac{\xi^2}{\sqrt{c^2t^2-\xi^2}} \; d\xi \\ 
& = 4\pi (ct)^2 \int_0^{r/(ct)} z^2 (1-z^2)^{-1/2} \; dz \\ 
& = 4\pi (ct)^2 \; \frac{1}{2} \left( \arcsin\left( \frac{r}{ct} \right) - \frac{r}{ct} \sqrt{1-\frac{r^2}{c^2t^2}} \right) \\
& = 2\pi (ct)^2 \arcsin\left( \frac{r}{ct} \right) - 2\pi r \sqrt{c^2t^2-r^2} ,
\endaligned  
\end{equation}
where we have used the easily checked equality: 
$$\int \frac{x^2}{\sqrt{1-x^2}} \; dx = \frac{1}{2} \left( \arcsin{x} - x\sqrt{1-x^2} \right) .$$
Then, by integrating the absolutely continuous part of (\ref{dens1}) over the ball $\in\mathcal B_r$ and taking into account (\ref{dens7}) and (\ref{dens8}). we have (for arbitrary $r<ct$):  
$$\aligned 
\text{Pr} \left\{ \bold X(t)\in\mathcal B_r \right\} & \sim e^{-\lambda t} \biggl[ \frac{\lambda}{4\pi c^2 t} \int\limits_{\mathcal B_r} \ln\left( \frac{ct+\Vert\bold x\Vert}{ct-\Vert\bold x\Vert} \right)  
\frac{d\bold x}{\Vert\bold x\Vert} + \frac{\lambda^2}{2\pi^2 c^2} \int\limits_{\mathcal B_r} \frac{d\bold x}{\sqrt{c^2t^2-\Vert\bold x\Vert^2}} + \frac{\lambda^3}{8\pi c^3} \int\limits_{\mathcal B_r} d\bold x \biggr] \\ 
& = e^{-\lambda t} \biggl[ \frac{\lambda}{4\pi c^2 t} \; 8\pi r ct \sum_{k=1}^{\infty} \frac{1}{4k^2-1} \; \left( \frac{r^2}{c^2t^2} \right)^k \\ 
& \qquad\qquad + \frac{\lambda^2}{2\pi^2 c^2} \biggl( 2\pi (ct)^2 \arcsin\left( \frac{r}{ct} \right) - 2\pi r \sqrt{c^2t^2-r^2} \biggr) + \frac{\lambda^3}{8\pi c^3} \; \frac{4}{3} \pi r^3 \biggr], 
\endaligned$$
and after some simple computations we finally arrive at the following asymptotic formula (for $r<ct$): 
\begin{equation}\label{dens9}
\aligned 
\text{Pr} \left\{ \bold X(t)\in\mathcal B_r \right\} & \sim e^{-\lambda t} \biggl[ \frac{2\lambda r}{c} \sum_{k=1}^{\infty} \frac{1}{4k^2-1} \; \left( \frac{r^2}{c^2t^2} \right)^k \\ 
& \qquad\qquad + \frac{\lambda^2 t^2}{\pi} \biggl( \arcsin\left( \frac{r}{ct} \right) - \frac{r}{ct} \sqrt{1 - \frac{r^2}{c^2t^2}} \biggr) + \frac{\lambda^3 r^3}{6 c^3} \biggr], \qquad t\to 0. 
\endaligned  
\end{equation}

\section{Estimate of the accuracy}

\numberwithin{equation}{section}

The error in asymptotic formula (\ref{dens1}) has the order $o(t^3)$. This means that, for small $t$, this formula yields a fairly 
good accuracy. To estimate it, let us integrate the function in square brackets of (\ref{dens1}) over the ball $\mathcal B_{ct}$. 

For the first term in square brackets of (\ref{dens1}) we have:  
\begin{equation}\label{est1}
\aligned 
\iiint\limits_{x_1^2+x_2^2+x_3^2\le c^2t^2} & \frac{\lambda}{4\pi c^2 t \Vert\bold x\Vert} \ln\left( \frac{ct+\Vert\bold x\Vert}{ct-\Vert\bold x\Vert} \right) \; dx_1 dx_2 dx_3 \\ 
& = \lambda t \; \iiint\limits_{x_1^2+x_2^2+x_3^2\le c^2t^2} \frac{1}{4\pi c^2 t^2 \Vert\bold x\Vert} \ln\left( \frac{ct+\Vert\bold x\Vert}{ct-\Vert\bold x\Vert} \right) \; dx_1 dx_2 dx_3 = \lambda t ,
\endaligned
\end{equation}
because the second integrand is the conditional density corresponding to the single change of direction (see \cite[the Theorem]{kol4} or \cite[formula (3.12)]{kol2}) and, therefore, the second integral is equal to 1.   

Applying \cite[formula 4.642]{gr}, we have for the second term in square brackets of (\ref{dens1}): 
\begin{equation}\label{est2}
\aligned 
\iiint\limits_{x_1^2+x_2^2+x_3^2\le c^2t^2} \frac{\lambda^2}{2\pi^2 c^2 \; \sqrt{c^2t^2-\Vert\bold x\Vert^2}} \; dx_1 dx_2 dx_3 & = \frac{\lambda^2}{2\pi^2 c^2} \; \iiint\limits_{x_1^2+x_2^2+x_3^2\le c^2t^2} \frac{dx_1 dx_2 dx_3}{\sqrt{c^2t^2-(x_1^2+x_2^2+x_3^2)}} \\ 
& = \frac{\lambda^2}{2\pi^2 c^2} \; \frac{2\pi^{3/2}}{\Gamma\left( \frac{3}{2} \right)} \; \int_0^{ct} \frac{\xi^2}{\sqrt{c^2t^2-\xi^2}} \; d\xi \\ 
& = \frac{2\lambda^2t^2}{\pi} \; \int_0^1 \frac{z^2}{\sqrt{1-z^2}} \; dz \\ 
& = \frac{\lambda^2t^2}{2} . 
\endaligned
\end{equation}

For the third term in square brackets of (\ref{dens1}) we get:
\begin{equation}\label{est3}
\iiint\limits_{x_1^2+x_2^2+x_3^2\le c^2t^2} \frac{\lambda^3}{8\pi c^3} \; dx_1 dx_2 dx_3 = \frac{\lambda^3}{8\pi c^3} \; \frac{4}{3} \; \pi c^3 t^3 = \frac{\lambda^3 t^3}{6} . 
\end{equation}

Hence, in view of (\ref{est1}), (\ref{est2}) and (\ref{est3}), the integral of the absolutely continuous part in asymptotic formula (\ref{dens1}) is: 
\begin{equation}\label{est4}
\aligned 
\tilde G(t) = \iiint\limits_{x_1^2+x_2^2+x_3^2\le c^2t^2} & e^{-\lambda t} \biggl[ \frac{\lambda}{4\pi c^2 t \Vert\bold x\Vert} \ln\left( \frac{ct+\Vert\bold x\Vert}{ct-\Vert\bold x\Vert} \right) + \frac{\lambda^2}{2\pi^2 c^2 \; \sqrt{c^2t^2-\Vert\bold x\Vert^2}} + \frac{\lambda^3}{8\pi c^3} \biggr] dx_1 dx_2 dx_3 \\ 
& = e^{-\lambda t} \left( \lambda t + \frac{\lambda^2t^2}{2} + \frac{\lambda^3 t^3}{6} \right) .
\endaligned
\end{equation}
Note that (\ref{est4}) can also be obtained by passing to the limit, as $r\to ct$, in asymptotic formula (\ref{dens9}).

On the other hand, according to (\ref{struc1}) and (\ref{densAC}), the integral of the absolutely continuous part of the transition density of the three-dimensional Markov random flight $\bold X(t)$ is 
\begin{equation}\label{est5}
G(t) = \int_{\mathcal B_{ct}} p^{(ac)}(\bold x,t) \; d\bold x = 1-e^{-\lambda t} . 
\end{equation}

The difference between the approximating function $\tilde G(t)$ and the exact function $G(t)$ given by (\ref{est4}) and (\ref{est5})   enables us to estimate the value of the probability generated by all the terms of the density aggregated in the term $o(t^3)$ of asymptotic relation (\ref{dens1}). 

The shapes of functions $G(t)$ and $\tilde G(t)$ on the time interval $t\in (0, \; 1)$ for the values of the intensity of switchings $\lambda =1, \; \lambda =1.5, \;\lambda =2, \;\lambda =2.5$ are presented in Figures 2 and 3. 

\begin{figure}[htbp]
\begin{minipage}[h]{0.45\linewidth}
\center{\includegraphics[width=1\linewidth]{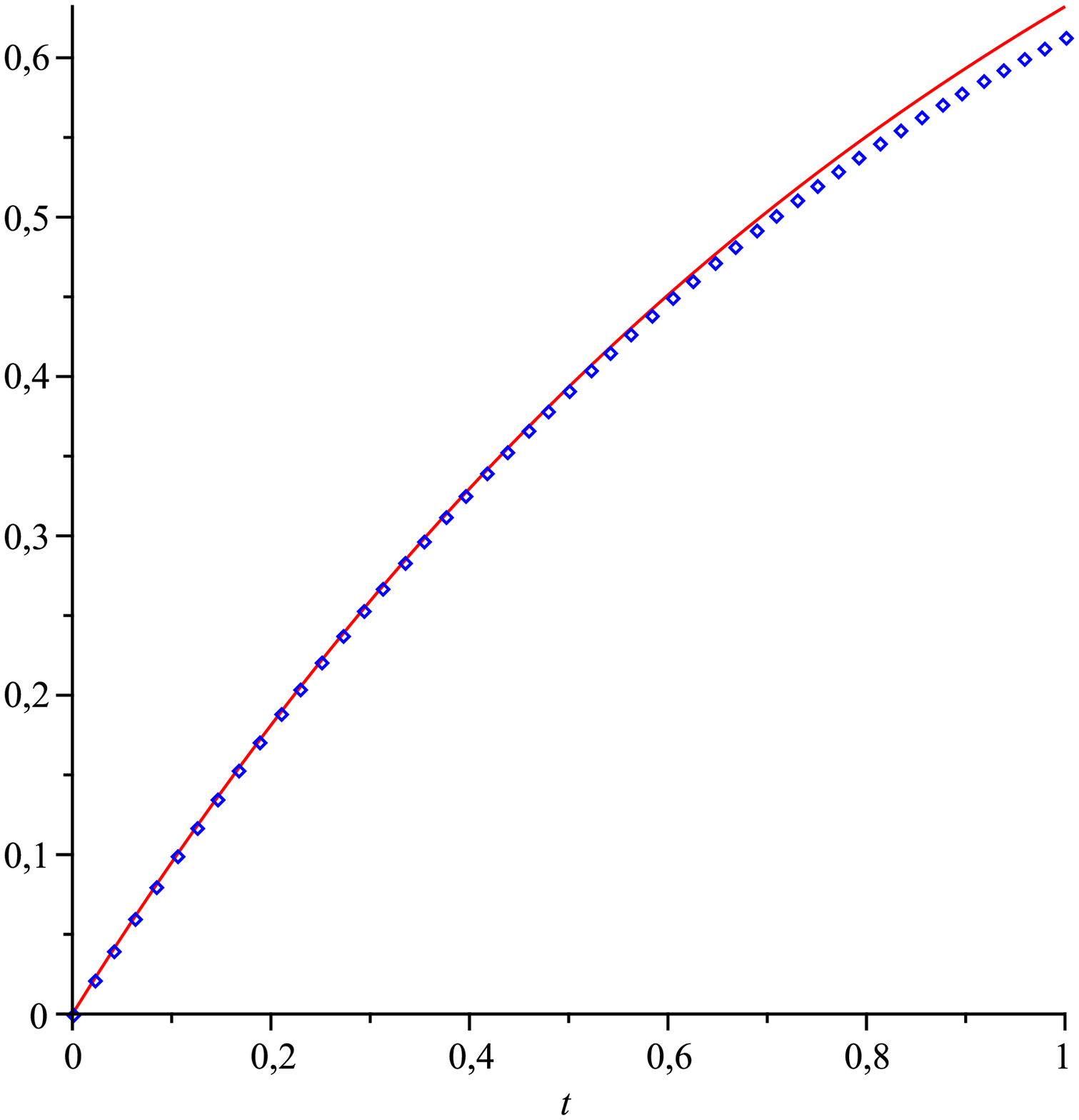}} \\
\end{minipage}
\hfill 
\begin{minipage}[h]{0.45\linewidth}
\center{\includegraphics[width=1\linewidth]{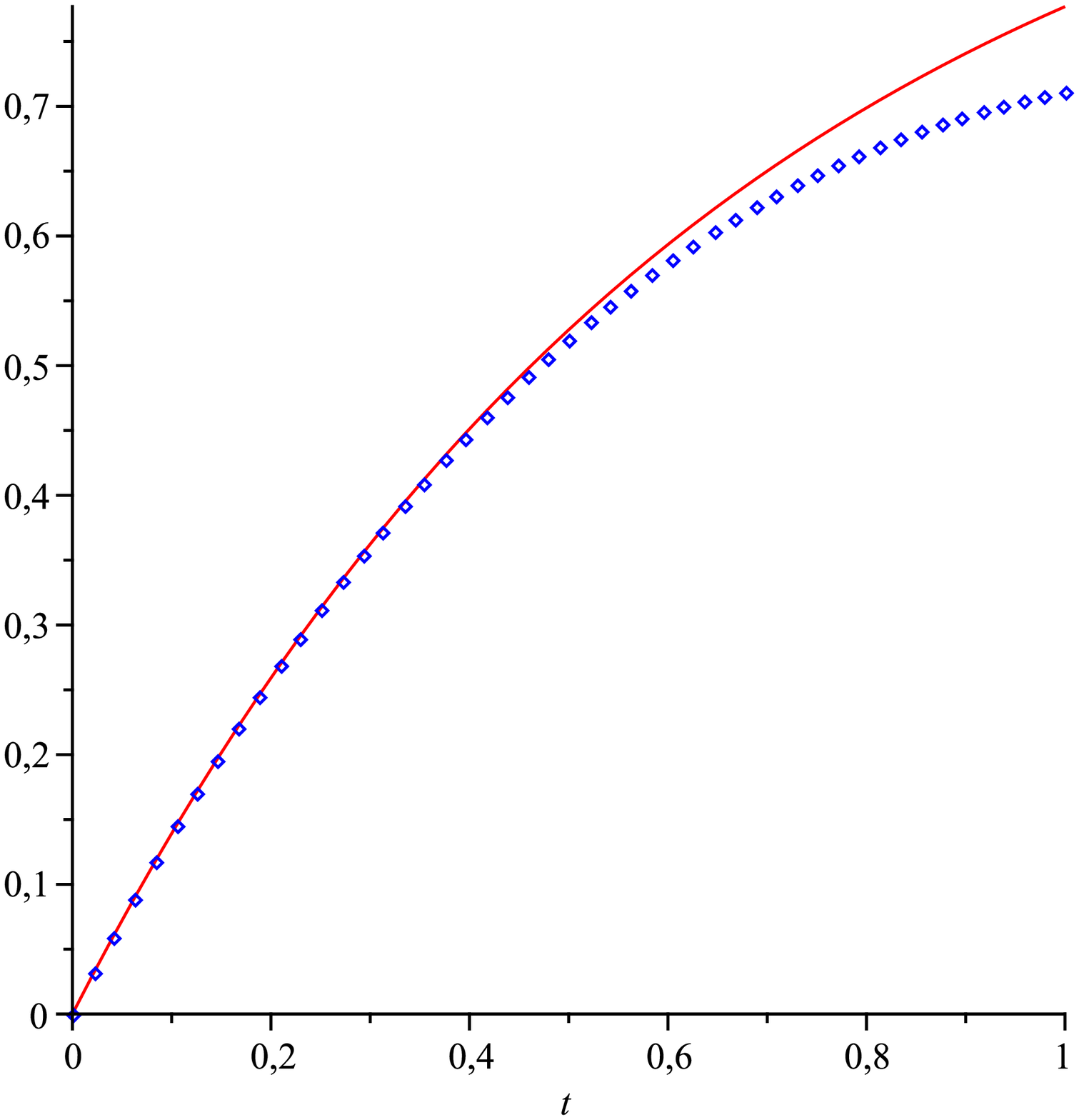}} \\ 
\end{minipage}
\caption{\it The shapes of functions $G(t)$ and $\tilde G(t)$ (point line) on the time interval $t\in (0, \; 1)$    
for the intensities $\lambda =1$ (left) and $\lambda = 1.5$ (right)}
\label{}
\end{figure}

\begin{figure}[htbp]
\begin{minipage}[h]{0.45\linewidth}
\center{\includegraphics[width=1\linewidth]{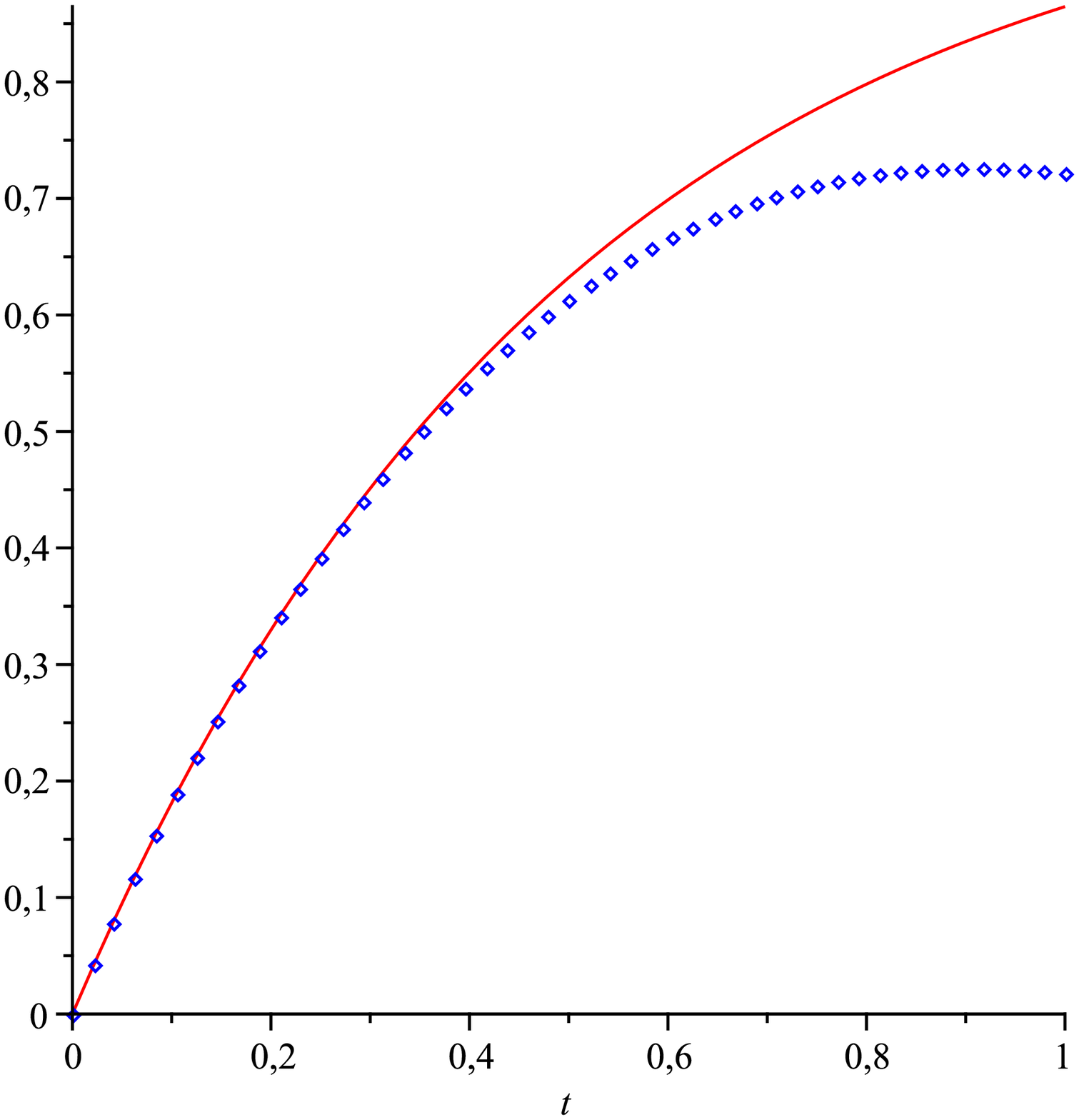}} \\
\end{minipage}
\hfill 
\begin{minipage}[h]{0.45\linewidth}
\center{\includegraphics[width=1\linewidth]{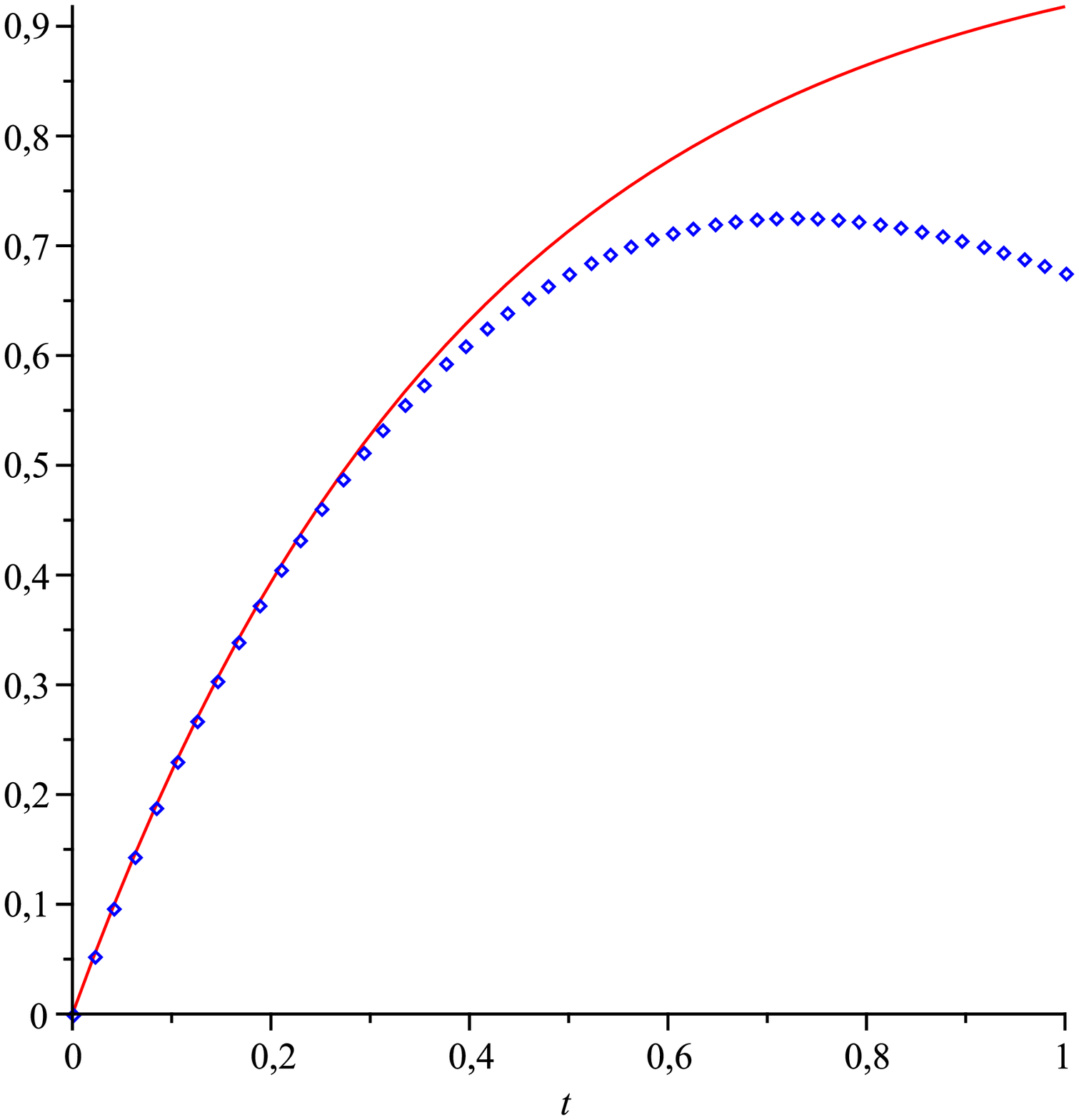}} \\ 
\end{minipage}
\caption{\it The shapes of functions $G(t)$ and $\tilde G(t)$ (point line) on the time interval $t\in (0, \; 1)$    
for the intensities $\lambda =2$ (left) and $\lambda = 2.5$ (right)}
\label{}
\end{figure}

We see that, for $\lambda =1$, the function $\tilde G(t)$ yields a very good coincidence with function $G(t)$ on the subinterval 
$t\in (0, \; 0.7)$ (Fig. 2 (left)), while for $\lambda =1.5$ (Fig. 2 (right)) such coincidence is good only on the subinterval 
$t\in (0, \; 0.5)$. The same phenomenon is also clearly seen in Figure 3 where, for $\lambda =2$, the function $\tilde G(t)$ 
yields a very good coincidence with function $G(t)$ on the subinterval $t\in (0, \; 0.4)$ (Fig. 3 (left)), while for $\lambda =2.5$ 
such good coincidence takes place only on the subinterval $t\in (0, \; 0.3)$ (Fig. 3 (right)). Thus, we can conclude that the greater is the intensity of switchings $\lambda$, the shorter is the subinterval of coincidence. This fact can easily be explained. Really, the greater is the intensity of switchings $\lambda$, the shorter is the time interval, on which no more than three changes of directions can occur with big probability. This means that, for increasing $\lambda$, the asymptotic formula (\ref{dens1}) yields a good accuracy on more and more small time intervals. However, for arbitrary fixed $\lambda$, there exists some $t_{\lambda}$ such that formula (\ref{dens1}) yields a good accuracy on the time interval $t\in (0, \; t_{\lambda})$ and the error of this approximation does not exceed $o(t_{\lambda}^3)$. This is the essence of the asymptotic formula (\ref{dens1}).

\bigskip

\begin{center}
{\Large\bf Appendices}
\end{center}

In the following appendices we establish some lemmas that have been used in the proofs of the above theorems. Note that some of them are of a separate mathematical interest because no similar results can be found in the mathematical handbooks. 

\begin{appendices}

\section{Auxiliary lemma}

{\bf Lemma A1.} {\it For arbitrary integer $n\ge 0$ and for arbitrary real $a \neq 0, -1, -2, \dots$, the following formula holds:}
\begin{equation}\label{appA1}
\sum_{k=0}^n \; \frac{ \Gamma\left( k+\frac{1}{2} \right) \; \Gamma\left( n-k+\frac{1}{2} \right)}{k! \; (n-k)! \; (2k+a)} 
= \frac{\pi \; \Gamma\left( \frac{a}{2} \right) \; \Gamma\left( n+\frac{a+1}{2} \right)}{(2n+a) \; \Gamma\left( \frac{a+1}{2} \right) \;
\Gamma\left( n+\frac{a}{2} \right)} , 
\end{equation}
$$n\ge 0, \;\; a \neq 0, -1, -2, \dots .$$

\vskip 0.2cm

{\it Proof.} Using the well-known relations for Pochhammer symbol 
\begin{equation}\label{apA3}
(-n)_k = \frac{(-1)^k  \; n!}{(n-k)!} , \;\; 0\le k\le n, \; n\ge 0, \qquad \frac{(x)_s}{(x+1)_s} = \frac{x}{x+s}, \;\; s>0,
\end{equation}
and the formula for Euler gamma-function  
\begin{equation}\label{appA2}  
\Gamma\left( k+\frac{1}{2} \right) = \frac{\sqrt{\pi}}{2^k} (2k-1)!!,  \qquad k=0,1,2,\dots, \quad (-1)!!=1,
\end{equation}
we can easily check that the sum on the left-hand side of (\ref{appA1}) is  
\begin{equation}\label{appA3}
\sum_{k=0}^n \; \frac{\Gamma\left( k+\frac{1}{2} \right) \; \Gamma\left( n-k+\frac{1}{2} \right)}{k! \; (n-k)! \; (2k+a)} = 
\frac{\sqrt{\pi} \; \Gamma\left( n+\frac{1}{2} \right)}{n! \; a} \;  _3F_2\left( -n, \frac{1}{2}, \frac{a}{2} ; -n+\frac{1}{2}, \frac{a}{2} +1; \; 1 \right) ,
\end{equation}
where 
$$_3F_2(\alpha_1,\alpha_2,\alpha_3; \; \beta_1,\beta_2; \; z) = \sum_{k=0}^{\infty} \frac{(\alpha_1)_k \; (\alpha_2)_k \; (\alpha_3)_k}{(\beta_1)_k \; (\beta_2)_k} \; \frac{z^k}{k!}$$
is the generalized hypergeometric function. According to \cite[item 7.4.4, page 539, formula 88]{pbm1} 
$$\aligned 
_3F_2\left( -n, \frac{1}{2}, \frac{a}{2} ; -n+\frac{1}{2}, \frac{a}{2} +1; \; 1 \right) & = \frac{\left( \frac{a+1}{2} \right)_n \; (1)_n}{\left( \frac{a}{2} +1 \right)_n \; \left( \frac{1}{2} \right)_n} \\
& = \frac{n! \; a \; \sqrt{\pi}}{\Gamma\left( n + \frac{1}{2} \right)} \; \frac{\Gamma\left( \frac{a}{2} \right) \; \Gamma\left( n + \frac{a+1}{2} \right)}{(2n+a) \; \Gamma\left( \frac{a+1}{2} \right) \; \Gamma\left( n + \frac{a}{2} \right)} .
\endaligned$$
Substituting this into (\ref{appA3}), we obtain (\ref{appA1}). The lemma is proved. $\square$

\section{Powers of the inverse tangent function}

In this appendix we derive series representations for some powers of the inverse tangent function that have been used in the proofs of the above theorems. Moreover, these results are of a more general mathematical interest because, to the best of the author's knowledge, there are no series representations, similar to (\ref{appB2}), (\ref{appB4}) and (\ref{appB6}) (see below), in mathematical handbooks, including \cite{gr}, \cite{korn}, \cite{pbm,pbm1}. 

\bigskip

{\bf Lemma B1.} {\it For arbitrary $z\in\Bbb C, \; |z|<\infty \; z\neq\pm i$, the following series representation holds:}
\begin{equation}\label{appB1}
\text{arctg} (z) = \frac{1}{\sqrt{\pi}} \; \frac{z}{\sqrt{1+z^2}} \; \sum_{k=0}^{\infty} \; \frac{ \Gamma\left( k+\frac{1}{2} \right) }{k! \; (2k+1)} \; \left( \frac{z^2}{1+z^2} \right)^k , \qquad |z|<\infty, \; z\neq\pm i . 
\end{equation}
{\it The series in} (\ref{appB1}) {\it is convergent uniformly in $z$.}

\vskip 0.2cm

{\it Proof.} Using well-known series representation of the inverse tangent function, see \cite[formula 1.644(1)]{gr}, as well as the formulas $(2k)!!= 2^k k! , \; k\ge 0,$ and (\ref{appA2}), we have (for $|z|<\infty, \; z\neq\pm i$):

$$\aligned
\text{arctg}(z) & = \frac{z}{\sqrt{1+z^2}} \; \sum_{k=0}^{\infty} \; \frac{(2k)!}{2^{2k} \; (k!)^2 \; (2k+1)} \; \left( \frac{z^2}{1+z^2} \right)^k \\ 
& = \frac{z}{\sqrt{1+z^2}} \; \sum_{k=0}^{\infty} \; \frac{(2k)!! \; (2k-1)!!}{ (2^k k!)^2 \; (2k+1)} \; \left( \frac{z^2}{1+z^2} \right)^k \\
& = \frac{z}{\sqrt{1+z^2}} \; \sum_{k=0}^{\infty} \; \frac{(2k-1)!!}{2^k \; k! \; (2k+1)} \; \left( \frac{z^2}{1+z^2} \right)^k \\
& = \frac{z}{\sqrt{1+z^2}} \; \sum_{k=0}^{\infty} \; \frac{\frac{\sqrt{\pi}}{2^k} (2k-1)!!}{2^k \; k! \; \frac{\sqrt{\pi}}{2^k} \; (2k+1)} \; \left( \frac{z^2}{1+z^2} \right)^k \\ 
& = \frac{1}{\sqrt{\pi}} \; \frac{z}{\sqrt{1+z^2}} \; \sum_{k=0}^{\infty} \; \frac{ \Gamma\left( k+\frac{1}{2} \right) }{k! \; (2k+1)} \; \left( \frac{z^2}{1+z^2} \right)^k ,
\endaligned$$
proving (\ref{appB1}). 

Since $\left\vert\frac{z^2}{1+z^2}\right\vert < 1$ for arbitrary $z\in\Bbb C, \; |z|<\infty, \; z\neq\pm i$, we get the inequality  
$$\left\vert\sum_{k=0}^{\infty} \; \frac{ \Gamma\left( k+\frac{1}{2} \right) }{k! \; (2k+1)} \; \left( \frac{z^2}{1+z^2} \right)^k \right\vert < \sum_{k=0}^{\infty} \; \frac{ \Gamma\left( k+\frac{1}{2} \right) }{k! \; (2k+1)} = \frac{\pi^{3/2}}{2},$$
proving the uniform convergence of the series in (\ref{appB1}). The lemma is proved. $\square$

\bigskip

{\bf Lemma B2.} {\it For arbitrary $z\in\Bbb C, \; |z|<\infty, \; z\neq\pm i$, the following series representation holds:}
\begin{equation}\label{appB2}
\bigl( \text{arctg}(z) \bigr)^2 = \frac{\sqrt{\pi}}{2} \left( \frac{z}{\sqrt{1+z^2}} \right)^2 \; \sum_{k=0}^{\infty} \; \frac{k!}{(k+1) \; \Gamma\left( k+\frac{3}{2} \right)} \; \left( \frac{z^2}{1+z^2} \right)^k , \qquad |z|<\infty, \; z\neq\pm i . 
\end{equation}
{\it The series in} (\ref{appB2}) {\it is convergent uniformly in $z$.}

\vskip 0.2cm

{\it Proof.} From (\ref{appB1}) it follows that 
\begin{equation}\label{appB3}
\bigl( \text{arctg}(z) \bigr)^2 = \frac{1}{\pi} \left( \frac{z}{\sqrt{1+z^2}} \right)^2 \; \sum_{k=0}^{\infty} \; \gamma_k \; \left( \frac{z^2}{1+z^2} \right)^k ,
\end{equation}
where the coefficients $\gamma_k$ are given by 
$$\gamma_k = \sum_{l=0}^k \frac{\Gamma\left( l+\frac{1}{2} \right) \; \Gamma\left( k-l+\frac{1}{2} \right) }{l! \; (k-l)! \; (2l+1) (2k-2l+1)} , \qquad k\ge 0 .$$
Since  
$$\frac{1}{(2l+1) (2k-2l+1)} = \frac{1}{2(k+1)} \left( \frac{1}{2l+1} + \frac{1}{2k-2l+1} \right) ,$$
then, taking into account the well-known formulas $z\Gamma(z) = \Gamma(z+1), \; \Gamma\left( \frac{1}{2} \right) = \sqrt{\pi}$, we have:

$$\aligned
\gamma_k & = \frac{1}{k+1} \sum_{l=0}^k \frac{\Gamma\left( l+\frac{1}{2} \right) \; \Gamma\left( k-l+\frac{1}{2} \right) }{l! \; (k-l)! \; (2l+1)} \quad (\text{see Lemma A1}) \\ 
& = \frac{1}{k+1} \; \frac{\pi \; \Gamma\left( \frac{1}{2} \right) \; \Gamma(k+1)}{(2k+1) \; \Gamma\left( k+\frac{1}{2} \right)} \\ 
& = \frac{\pi^{3/2} \; k!}{ 2(k+1) \; \left( k+\frac{1}{2} \right) \Gamma\left( k+\frac{1}{2} \right)} \\
& = \frac{\pi^{3/2} \; k!}{2(k+1) \; \Gamma\left( k+\frac{3}{2} \right)} . 
\endaligned$$
Substituting these coefficients into (\ref{appB3}) we obtain (\ref{appB2}). The uniform convergence of the series in formula (\ref{appB2}) can be established similarly to that of Lemma B1. This completes the proof of the lemma. $\square$

\bigskip

{\bf Lemma B3.} {\it For arbitrary $z\in\Bbb C, \; |z|<\infty, z\neq\pm i$, the following series representation holds:}
\begin{equation}\label{appB4}
\aligned 
\bigl( \text{arctg}(z) \bigr)^3 & = \frac{1}{\sqrt{\pi}} \left( \frac{z}{\sqrt{1+z^2}} \right)^3 \\ 
& \times \sum_{k=0}^{\infty} \frac{\Gamma\left( k+\frac{1}{2} \right)}{k! \; 
(2k+1)} \; _5F_4\left( 1,1,1,-k,-k-\frac{1}{2}; \; -k+\frac{1}{2}, -k+\frac{1}{2}, \frac{3}{2}, 2; \; 1 \right) \left( \frac{z^2}{1+z^2} \right)^k ,  
\endaligned
\end{equation}
{\it where}
\begin{equation}\label{hypergeom54}
_5F_4(a_1,a_2,a_3,a_4,a_5; \; b_1,b_2,b_3,b_4; \; z) = \sum_{k=0}^{\infty} \frac{(a_1)_k \; (a_2)_k \; (a_3)_k \; (a_4)_k \; (a_5)_k}{(b_1)_k \; (b_2)_k \;(b_3)_k \;(b_4)_k} \; \frac{z^k}{k!}
\end{equation} 
{\it is the generalized hypergeometric function. The series in} (\ref{appB4}) {\it is convergent uniformly in $z$.}
\vskip 0.2cm

{\it Proof.} From (\ref{appB1}) and (\ref{appB2}) it follows that  
\begin{equation}\label{appB5}
\bigl( \text{arctg}(z) \bigr)^3 = \frac{1}{2} \left( \frac{z}{\sqrt{1+z^2}} \right)^3 \; \sum_{k=0}^{\infty} \; \gamma_k \; \left( \frac{z^2}{1+z^2} \right)^k ,
\end{equation}
where the coefficients $\gamma_k$ are given by 
$$\gamma_k = \sum_{l=0}^k \frac{l! \;\; \Gamma\left( k-l+\frac{1}{2} \right)}{(l+1) \; (k-l)! \; (2k-2l+1) \; \Gamma\left( l+\frac{3}{2} \right)} , \qquad k\ge 0 .$$
Applying (\ref{apA3}), (\ref{appA2}) and the formula 
$$\Gamma\left( \frac{1}{2} - k \right) = \frac{(-1)^k \; \sqrt{\pi} \; 2^k}{(2k-1)!!} , \qquad k\ge0,$$
after some simple computations, we arrive at the relation   
$$\gamma_k = \frac{2\Gamma\left( k+\frac{1}{2} \right)}{k! \; \sqrt{\pi} (2k+1)} \; _5F_4\left( 1,1,1,-k,-k-\frac{1}{2}; \; -k+\frac{1}{2}, -k+\frac{1}{2}, 
\frac{3}{2}, 2; \; 1 \right) , \qquad k\ge 0.$$
Substituting these coefficients into (\ref{appB5}) we obtain (\ref{appB4}). The lemma is proved. $\square$

\bigskip

{\bf Lemma B4.} {\it For arbitrary $z\in\Bbb C, \; |z|<\infty, \; z\neq\pm i$, the following series representation holds:}
\begin{equation}\label{appB6}
\bigl( \text{arctg}(z) \bigr)^4 = \frac{\pi}{2} \left( \frac{z}{\sqrt{1+z^2}} \right)^4 \; \sum_{k=0}^{\infty} \gamma_k \left( \frac{z^2}{1+z^2} \right)^k , \qquad |z|<\infty, \; z\neq\pm i , 
\end{equation}
{\it where the coefficients $\gamma_k$ are given by the formula}
$$\gamma_k = \frac{1}{k+2} \sum_{l=0}^k \frac{l! \; (k-l)!}{(l+1) \; \Gamma\left( l+\frac{3}{2} \right) \; \Gamma\left( k-l+\frac{3}{2} \right)} , \qquad k\ge 0.$$
{\it The series in} (\ref{appB6}) {\it is convergent uniformly in $z$.}

\vskip 0.2cm

{\it Proof.} According to Lemma B2, we have: 
\begin{equation}\label{appB7}
\bigl( \text{arctg}(z) \bigr)^4 = \frac{\pi}{4} \left( \frac{z}{\sqrt{1+z^2}} \right)^4 \; \sum_{k=0}^{\infty} \xi_k \left( \frac{z^2}{1+z^2} \right)^k ,
\end{equation}
where the coefficients $\xi_k$ are: 
$$\aligned
\xi_k & = \sum_{l=0}^k \frac{l! \; (k-l)!}{(l+1) (k-l+1) \; \Gamma\left( l+\frac{3}{2} \right) \; \Gamma\left( k-l+\frac{3}{2} \right)} \\
& = \frac{1}{k+2} \sum_{l=0}^k \frac{l! \; (k-l)!}{\Gamma\left( l+\frac{3}{2} \right) \; \Gamma\left( k-l+\frac{3}{2} \right)} \left[ \frac{1}{l+1} + \frac{1}{k-l+1} \right]\\
& = \frac{2}{k+2} \sum_{l=0}^k \frac{l! \; (k-l)!}{(l+1) \; \Gamma\left( l+\frac{3}{2} \right) \; \Gamma\left( k-l+\frac{3}{2} \right)} . 
\endaligned$$
Substituting this into (\ref{appB7}), we get the statement of the lemma. $\square$

\end{appendices}

\end{document}